\documentclass{tac}
\usepackage{amsmath}
\usepackage{amssymb}
\usepackage{tikz-cd}

\usepackage[colorlinks=true]{hyperref}
\hypersetup{allcolors=[rgb]{0.1,0.1,0.4}}

\let\pf\proof
\let\epf\endproof

\mathrmdef[relfin]{RelFin}
\mathrmdef[finset]{FinSet}
\mathrmdef[grpd]{Grpd}
\mathrmdef[cat]{Cat}
\mathrmdef[sym]{Sym}
\mathrmdef[aut]{Aut}

\mathrmdef{Hom}
\mathrmdef{hom}
\mathrmdef{inj}
\mathrmdef{core}
\mathrmdef{im}
\mathrmdef{exp}
\mathrmdef{id}

\newcommand{\Fq}{\mathbb F_q}

\renewcommand{\cal}[1]{\mathcal{#1}}

\newtheorem{conj}{Conjecture}

\title{Combinatorics in (2,1)-categories}
\author{Krista Zehr}
\date{February 2026}
\copyrightyear{2026}
\thanks{Thanks to Justin Curry and Marco Varisco for their feedback and 
for suggesting several of the questions investigated in this project.}

\begin{document}

\maketitle
\begin{abstract}
	Groupoid cardinality is an invariant of locally finite groupoids 
	which has many of the properties of the cardinality of finite sets, 
	but which takes values in all non-negative real numbers, and 
	accounts for the morphisms of a groupoid. Several results on 
	groupoid cardinality are proved, analogous to the relationship 
	between cardinality of finite sets and i.e. injective or surjective 
	functions. We also generalize to a broad class of (2,1)-categories 
	a famous theorem of Lov\'asz which characterizes the isomorphism 
	type of relational structures by counting the number of 
	homomorphisms into them.
\end{abstract}

\section{Introduction}
When studying mathematical objects defined as finite sets with extra 
structure, cardinality is a very powerful tool. While it is a complete 
invariant of finite sets up to bijection, this is not the main source 
of its power. For example, monomorphisms or epimorphisms enforce an 
order on the cardinalities of their source and target objects. 
Moreover, the cardinalities of the source and target of a morphism can 
enforce restrictions on the properties of that morphism. Categorical 
constructions such as products or coproducts reduce to algebraic 
equations involving the cardinalities of the component objects.

It was shown in \cite{lovasz1} that in a category of finite relational 
structures, it is enough to count the number of homomorphisms into an 
object from every other object to determine its isomorphism type. This 
result can be extended to many other locally finite categories. One way 
of understanding this result is via the Yoneda embedding of a category 
$\cal C \to [\cal C^{op}, \finset]$. That is, given two representable 
functors $h_A$ and $h_B$, and a \emph{natural} isomorphism between 
them, the Yoneda lemma states that $A$ and $B$ must be isomorphic in 
$\cal C$. However, if we weaken the natural isomorphism to just a 
family of isomorphisms $\Hom(X,A) \cong \Hom(X,B)$ for each object $X$, 
then this is no longer true. Since cardinality is a complete invariant 
of sets up to isomorphism, Lov\'asz's result and its generalizations 
gives conditions where we can drop the naturality assumption and still 
conclude that $A \cong B$.

Lov\'asz originally used this lemma to show that products of finite 
relational structures obey a cancellation property. That is, if 
$A\times B \cong A \times C$, then $B \cong C$, which follows directly 
from the universal property of products. Beyond this result there have 
been other applications of counting homomorphisms in this way. For 
example, the Weisfeiler-Leman isomorphism test is a polynomial time 
algorithm that can distinguish between isomorphism classes of finite 
graphs. It was shown in \cite{dgr} that this test fails to distinguish 
between graphs $G$ and $G'$ if and only if the number of homomorphisms 
$T \to G$ is the same as the number of homomorphisms $T \to G'$ for all 
\emph{trees} $T$. Moreover, this generalizes to the $k$-dimensional 
version of the Weisfeiler-Leman test by counting homomorphisms 
from graphs $\Gamma$ of treewidth at most $k$.

When discussing higher categories, the morphisms between any two 
objects no longer form a set, but a category, so an alternative 
invariant should be used. One invariant is the Euler characteristic 
of a category (c.f. \cite{leinster}). The Euler characteristic is an 
invariant up to equivalence (in fact, up to adjunction), and satisfies 
algebraic formulas for products and coproducts identical to those for 
cardinality of finite sets. The downsides are that it is hard to 
compute in general, and the categories that have well-defined Euler 
characteristic are not easy to classify. However, if the category in 
question is a groupoid, then the Euler characteristic is known as 
groupoid cardinality, and has a simple formula.

Groupoid cardinality essentially counts objects, but only up to 
isomorphism. So the cardinality of a discrete groupoid is equal to the 
number of objects, and the cardinality of an equivalence relation is 
equal to the number of equivalence classes. Moreover, when the groupoid 
is not a preorder, we may have fractional cardinality. For instance, if 
$G$ is a finite group, then the cardinality of the delooping groupoid 
$BG$ is the reciprocal of the order of $G$. In this way, groupoid 
cardinality can be a useful invariant when we want to consider 
symmetries or equivalences between objects.

One of the goals of this paper is to further study the properties of 
groupoid cardinality analogous to the useful properties of set 
cardinality that were previously mentioned. In Section 
\ref{sec_funct_grpd_card}, we study the cardinalities of exponential 
objects in the category of groupoids $\grpd$. These cardinalities are 
hard to compute in general, and have much more complex formulas than 
the simple exponential formula for the cardinality of the set of 
functions between finite sets. In Section \ref{sec_morphisms_tame}, we 
relate properties of functors to groupoid cardinality.

Moving beyond the the category of groupoids, we investigate other 
(2,1)-categories. In such categories, the morphisms between objects 
form a groupoid, so we may calculate their groupoid cardinality. One 
example is the category of stuff types (see \cite{morton}), a 
generalization of Joyal's combinatorial species, which are themselves a 
sort of categorification of generating functions in combinatorics. 
Stuff types are simply functors from an arbitrary groupoid into the 
groupoid of finite sets and bijections. This picture can be generalized 
very broadly, replacing the groupoid of finite sets with any locally 
finite groupoid. One specific example mentioned in Section 
\ref{sec_funct_grpd_card} is the groupoid of finite dimensional vector 
spaces over a finite field and the invertible linear maps between them. 
Section \ref{sec_rel_fin} discusses such (2,1)-categories as well as a 
certain finiteness condition on the functors in them that is essential 
for combinatorial analysis.

The main theorem of this paper is a generalization of Lov\'asz's 
theorem to such categories of functors, and is proven in Section 
\ref{sec_comb_cat}. Finally, in Section \ref{sec_homotopy} we discuss 
possible future directions of this work to ($\infty$,1)-categories, 
using the notion of homotopy cardinality of an $\infty$-groupoid, which 
was first introduced in \cite{bd}.

\section{Background}

\subsubsection{Groupoid cardinality}
Groupoid cardinality was first introduced in \cite{bd} as a generalization of cardinality of finite sets which can take values of any non-negative real number.

Let $\cal G$ be an essentially small groupoid. Denote the set of all isomorphism classes of $\cal G$ by $\overline{\cal G}$. If $x$ is an object of $\cal G$, then $[x]$ denotes the isomorphism class of $x$, and $\cal G_x$ is the group of automorphisms of $x$. For the remainder of this paper, we will denote the cardinality of a finite set $X$ by $\#X$, in order to avoid confusion with groupoid cardinality.

\definition
    Let $\cal G$ be a locally finite, essentially small groupoid. Then 
    the \textbf{groupoid cardinality} of $\cal G$ is
    \[
    |\cal G| = \sum_{[x] \in \overline{\cal G}} \frac1{\#\cal G_x}
    \]
    We call such a groupoid \textbf{tame} if $|\cal G| < \infty$.
\enddefinition
The following properties of groupoid cardinality are easy to verify.
\proposition \label{card_facts}
    Let $\mathcal {G,H}$ be tame groupoids.
    \begin{enumerate}
        \item [(a)] If $\cal G \simeq \cal H$ then $|\cal G| = |\cal H|$
        \item [(b)] $|\cal G \times \cal H| = |\cal G||\cal H|$
        \item [(c)] $|\cal G \sqcup \cal H| = |\cal G| + |\cal H|$
        \item [(d)] $|\cal G| = 0$ if and only if $\cal G$ is the empty groupoid.
        \item [(e)] If $\cal G$ is finite, then $|\cal G|$ is rational.
    \end{enumerate}
\endproposition

\proposition
Let $G$ be a finite group, acting on a finite set $X$. Then the cardinality of the action groupoid $X // G$ is
\[
|X // G| = \frac{\#X}{\#G}
\]
\endproposition

Clearly, groupoid cardinality has similar formal properties to the cardinality of finite sets. Note that the category of functors between groupoids is again a groupoid, so it is natural to ask if the equation
\[
|\mathcal{G^H}| = |\cal G|^{|\cal H|}.
\]
holds for tame groupoids, as it does for the cardinality of finite 
sets. However, this fails miserably. In general the functor category 
between two tame groupoids is not even locally finite.

When $\cal H$ is finite, the functor category is locally finite, 
but this equation still fails to hold. For example, consider the 
functors from the cyclic group $C_2$ to the discrete groupoid with 2 
objects. This is obviously a finite groupoid. However, if the above 
equation were true, then the cardinality would be $\sqrt 2$ which is 
not rational, contradicting Proposition~\ref{card_facts}(e). We will 
investigate cardinalities of functor groupoids in Section 
\ref{sec_funct_grpd_card}.

\subsection{Stuff Types}
\definition
A \textbf{stuff type} is a pair $(\cal G, \Phi)$ where $\cal G$ is a 
groupoid and $\Phi: \cal G \to \core(\finset)$ is a functor.
\enddefinition
For the sake of brevity, we will denote $\cal B := \core(\finset)$ in 
this section.

Stuff types were first introduced in \cite{bd}, as a categorification 
of generating functions. The objects of $\cal G$ can be interpreted as 
finite sets with some extra ``stuff" attached, and the morphisms can be 
interpreted as bijections that preserve the attached ``stuff". The 
functor $\Phi$ can be interpreted as a forgetful functor sending 
objects to their underlying sets, and morphisms to their underlying 
bijection.

\definition
The \textbf{generating function} of a stuff type $\Phi: \cal G \to 
\cal B$ is the formal power series
\[
|\Phi|(z) := \sum_{n \in \mathbb N} |\Phi^{-1}(n)|z^n,
\]
where $\Phi^{-1}(n))$ denotes the full subgroupoid of $\cal G$ on all 
objects that get sent to a set of cardinality $n$ by $\Phi$.
\enddefinition

\definition
The \textbf{sum} of two stuff types $\Phi:\cal G \to \cal B$ and $\Psi: 
\cal H \to \cal B$ is the functor $\Phi + \Psi : \cal G \sqcup \cal H 
\to \cal B$, which is equal to the coproduct of $\Phi$ and $\Psi$ in 
the slice (2,1)-category $\grpd/\cal B$. 
\enddefinition

\definition
The \textbf{product} (often referred to as the Cauchy product to 
distinguish it from the categorical product) of two stuff types 
$\Phi:\cal G \to \cal B$ and $\Psi: \cal H \to \cal B$ is the functor 
$\Phi\Psi: \cal G \times \cal H \to \cal B$ which acts on objects and 
morphisms by $(g,h) \mapsto \Phi(g)\sqcup \Psi(h)$
\enddefinition

These two definitions are analogous to the sum and product operations 
on generating functions. That is, suppose we have stuff types $\Phi: 
\cal G \to \cal B$ and $\Psi: \cal H \to \cal B$. An object of the 
domain of the sum $\Phi + \Psi : \cal G \sqcup \cal H \to \cal B$ is a 
finite set equipped with either $\Phi$-stuff or $\Psi$-stuff. An object 
in the domain of the product $\Phi\Psi :\cal G \times \cal H \to \cal 
B$ is a finite set partitioned into two subsets, one equipped with 
$\Phi$-stuff, and the other equipped with $\Psi$-stuff. Moreover, sums 
and products behave nicely when taking generating functions, as shown 
in the following proposition:

\proposition
Let $\Phi:\cal G \to \cal B$ and $\Psi: \cal H \to \cal B$ be stuff 
types.
\begin{enumerate}
	\item $|\Phi + \Psi|(z) = |\Phi|(z) + |\Psi|(z)$
	\item $|\Phi\Psi|(z) = |\Phi|(z) \cdot |\Psi|(z)$
\end{enumerate}
\endproposition

\pf
\begin{enumerate}
	\item Follows trivially from the definition and 
	Proposition~\ref{card_facts}(c).
	\item The coefficient of $z^n$ in $|\Phi\Psi|(z)$ is the 
	cardinality of $(\Phi\Psi)^{-1}(n)$. This groupoid is the full 
	subgroupoid of $\cal G \times \cal H$ on all objects $(X,Y)$ such 
	that $\Phi(X) \sqcup \Psi(Y)$ is an $n$-element set. This means 
	that if we let $i := \#\Phi(X)$ and $j := \#\Psi(Y)$, then $i+j = 
	n$, so 
	\[
		(\Phi\Psi)^{-1}(n) = \bigsqcup_{i+j = n} \Phi^{-1}(i)\times 
		\Psi^{-1}(j).
	\]
	Taking the cardinality, and applying Proposition~\ref{card_facts}, 
	we get
	\[
	|(\Phi\Psi)^{-1}(n)| = \sum_{i+j=n} |\Phi^{-1}(i)||\Psi^{-1}(j)|
	\]
	Which is exactly the formula for the $n$th coefficient of the 
	product of formal power series $|\Phi|(z)|\Psi|(z)$.
\end{enumerate}
\epf

Note that the previous definitions and results do not rely on any facts 
about $\core(\finset)$ other than the fact that the isomorphism 
classes can be identified with natural numbers in such away that 
coproducts correspond to addition under this identification. Therefore, 
we could replace $\core(\finset)$ with another groupoid with 
the same property, such as finite dimensional vector spaces over a 
field $k$.

\subsection{Lov\'asz's homomorphism counting lemma}
In this section we will discuss the result in \cite{lovasz1} that 
distinguishes isomorphism types by counting homomorphisms into an 
object. The proof of this statement is elementary, and very simple to 
understand, so we will provide it in full here.

\definition
A \textbf{finite relational structure} is a tuple $(A, R_1, \dots, 
R_n)$ where $A$ is a finite set and $R_i \subseteq A^{a_i}$ are 
relations on $A$ with arity $a_i$ (we assume $n$ and the arities $a_i$ 
to be fixed). A \textbf{homomorphism of relational structures} from 
$(A,R)$ to $(B, S)$ is a function $f:A \to B$ between the underlying 
sets such that $f(R_i) \subseteq S_i$ for all $i$. That is, $f$ 
preserves the relations.
\enddefinition

\theorem \label{lovasz_lemma} \cite{lovasz1}
Let $(A,R)$, $(B,S)$ be finite relational structures, and let 
$\hom((C,T),(A,R))$ denote the number of homomorphisms from 
$(C,T)$ to $(A,R)$. If for every structure $(C,T)$, we have 
$\hom((C,T),(A,R)) = \hom((C,T),(B,S))$, then $(A,R) 
\cong (B,S)$.
\endtheorem

\pf
In the first half of the proof, we prove by induction on the 
cardinality of $C$ that if $\hom((C,T),(A,R)) = \hom((C,T),(B,S))$ for 
all $(C,T)$, then it is also true that $\inj((C,T),(A,R)) = 
\inj((C,T),(B,S))$, where $\inj((C,T),(A,R))$ is the number of 
injective homomorphisms from $(C,T)$ to $(A,R)$.

For the base case, if $C = \varnothing$, then $\inj((C,T),(A,R)) = 
\inj((C,T),(B,S)) = 1$ trivially.

For the inductive case, assume that $\inj((C,T),(A,R)) = 
\inj((C,T),(B,S))$ for all $(C,T)$ such that $|C| < N$ for some 
$N > 0$. Consider that every homomorphism of relational structures 
$f:(C,T) \to (A,R)$ decomposes uniquely as $f = m\circ e$ where $m$ is 
an injective homomorphism and $e$ is a surjective homomorphism. (This 
is easy to check). Therefore, we have the equation
\[
\hom((C,T),(A,R)) = \sum_{\theta} \inj((C/\theta,T),(A,R)),
\]
where the sum is over equivalence relations $\theta$ on $C$, 
interpreting $C/\theta$ as the image of a homomorphism in $A$.
Now, fix $(C,T)$ with $|C| = N$. Then we have
	$$\mld 0 &= \hom((C,T),(A,R)) - \hom((C,T),(B,S))\\ = \sum_{\theta} 
	\inj((C/\theta,T),(A,R)) - \inj((C/\theta,T),(B,S))\\ = 
	\inj((C,T),(A,R)) - \inj((C,T),(B,S)) + \sum_{\theta \text{ 
	nontrivial}} \inj((C/\theta,T),(A,R)) - \inj((C/\theta,T),(B,S))$$
By our induction hypothesis, the terms in the sum on the right are all 
zero, therefore $\inj((C,T),(A,R)) = \inj((C,T),(B,S))$.

For the second half of the proof, we set $(C,T) = (A,R)$, so that 
$\inj((A,R),(A,R)) = \inj((A,R),(B,S)) \neq 0$. So 
there exists some injective homomorphism from $(A,R)$ to $(B,S)$. 
Setting $(C,T) = (B,S)$, we can see that there also exists an injective 
homomorphism from $(B,S)$ to $(A,R)$. Since these homomorphisms are 
functions of sets, we can see that they must be isomorphisms, which 
proves that $(A,R) \cong (B,S)$.
\epf

This theorem is very general, and yet the proof is very easy. It is 
easy to see how the basic structure of the proof could be generalized 
to other areas, and indeed, a more categorical description was 
developed in \cite{lovasz2} and \cite{pultr}.

\subsection{Factorization systems in (2,1)-categories} 
\label{factorization}
The proof of Theorem \ref{lovasz_lemma} relies on the existence of image factorizations. In this section we will review the theory of factorization systems in a 2-category. The following definition is from \cite{kv00}.
\definition \label{def_fact_syst}
    A \textbf{factorization system} in a (2,1)-category $\cal C$ is a pair $(\cal E, \cal M)$ of classes of 1-morphisms such that:
    \begin{enumerate}
        \item $\cal E$ and $\cal M$ contain all equivalences and are closed under composition with equivalences,
        \item $\cal E$ and $\cal M$ are closed under 2-isomorphism classes of 1-morphisms,
        \item For every 1-morphism $f:X\to Y$, there exist $e$ in $\cal 
        E$, $m$ in $\cal M$ and a 2-cell $\eta: me \Rightarrow f$.
        \item Finally, $\cal E, \cal M$ are required to satisfy an 
        orthogonality or "fill-in" condition. Given any diagram:
        \[\begin{tikzcd}
        	A & B \\
        	C & D
        	\arrow["e", from=1-1, to=1-2]
        	\arrow["u"', from=1-1, to=2-1]
        	\arrow["m"', from=2-1, to=2-2]
        	\arrow["v", from=1-2, to=2-2]
        	\arrow["\varphi"', Rightarrow, from=1-2, to=2-1]
        \end{tikzcd}\]
        with $e$ in $\cal E$ and $m$ in $\cal M$, there exists $w : B \to C$ and 2-morphisms $\alpha: we \Rightarrow u$ and $\beta: mw \Rightarrow v$ such that $\varphi\circ (\beta * 1_e) = 1_m * \alpha$. Moreover, if $(w', \alpha', \beta')$ is another such fill-in, then then there exists a unique $\psi: w \Rightarrow w'$ such that $\alpha'(\psi * 1_e) = \alpha$ and $\beta'(1_m * \psi) = \beta$.
    \end{enumerate}
\enddefinition
By Proposition 9.3 in \cite{kv00}, we can consider $\cal E$ and $\cal M$ as subcategories of the underlying 1-category of $\cal C$. Furthermore, factorizations of morphisms are unique in the following sense:
\proposition
    If $(e:X\to U,m:U \to Y,\eta)$ and $(e':X\to U', m':U'\to Y,\eta')$ 
    are two factorizations of $f$, then there exits an equivalence $w:U 
    \to U'$ and 2-cells $\varphi : we \Rightarrow e', \psi : m'w 
    \Rightarrow m$ such that $\eta'(1_{m'}\ast\varphi) = 
    \eta(\psi\ast1_e)$. Moreover $w$ is unique up to a unique coherent 
    2-isomorphism.
\endproposition

The proof of this proposition is a straightforward exercise in using 
the orthogonality condition.

This definition of a factorization system is reminiscent of epi-mono 
factorization systems in ordinary 1-categories. Indeed, it makes sense 
to consider $\cal E$ to be akin to quotient maps, and $\cal M$ to be 
akin to subobject inclusions. Often, quotient maps and subobject 
inclusions have a variety of nice properties beyond being factors of 
morphisms, which motivates the following definition:

\definition[{\cite[Definition 3.1]{dv03}}]\label{fully_co_faith}
Let $f: C \to C'$ be a morphism in a 2-category $\cal C$
\begin{enumerate}
    \item We say $f$ is \textbf{fully faithful} if for all objects $X$ of $\cal C$, the induced functor $f_*:\cal C(X,C) \to \cal C(X,C')$ is fully faithful.
    \item We say $f$ is \textbf{fully cofaithful} if for all objects $X$ of $\cal C$, the induced functor $f^*:\cal C(C',X) \to \cal C(C,X)$ is fully faithful.
\end{enumerate}
\enddefinition

In section \ref{sec_morphisms_tame}, we provide a number of results 
relating groupoid cardinality to properties of functors between tame 
groupoids (i.e. full, faithful, essentially surjective). Since Theorem 
\ref{main_thm} refers to more general (2,1)-categories than $\grpd$, 
these properties of morphisms in a 2-category will allow us to use 
these results in this more general setting.

For small categories (and therefore, small groupoids) these definitions agree with well-understood properties of functors.
\proposition[{\cite[Proposition 7.8]{dv03}}]\label{dv_prop_proper_fact}
Let $F:\cal C \to \cal D$ be a functor between small categories. Then
\begin{enumerate}
    \item $F$ is fully faithful in the sense of Definition~\ref{fully_co_faith} if and only if $F$ is fully faithful as a functor.
    \item $F$ is fully cofaithful if $F$ is essentially surjective and full as a functor.
\end{enumerate}
\endproposition

Factorizing morphisms into a composite of two morphisms, one a "quotient" and the other an "embedding", is a common and well-understood theme in mathematics. When moving to 2-categories, an extra layer of subtlety arises. Namely, there is more than one obvious choice of factorization system, leading to the following definition:

\definition
    A \textbf{ternary factorization system} is a pair of factorization systems $(\cal E, \widetilde{\cal M})$ and $(\widetilde{\cal E}, \cal M)$ such that $\cal E \subseteq \widetilde{\cal E}$ and $\cal M \subseteq \widetilde{\cal M}$.
\enddefinition
We call this a ternary factorization system because each morphism 
factors uniquely as three morphisms rather than two. Given a morphism 
$\varphi$, we first factor $\varphi$ as $\tilde{m}e$ with $e$ in $\cal 
E$ and $\tilde{m}$ in $\widetilde{\cal M}$, and then factor $\tilde{m}$ 
as $mf$ for $m$ in $\cal M$ and $f$ in $\widetilde{\cal E}$. 
Alternatively, we may factor $\varphi$ with $(\widetilde{\cal E},\cal 
M)$ first and then with $(\cal E,\widetilde{\cal M})$. The 
orthogonality condition ensures that both of these factorizations are 
equivalent. Thus, we may consider a ternary factorization system to be 
a triple of subcategories $(\cal E, \cal F, \cal M)$ with $\cal F = 
\widetilde{\cal E} \cap \widetilde{\cal M}$. With this definition, each 
morphism factors as a composition $mfe$ with $e$ in $\cal E$, $f$ in 
$\cal F$, and $m$ in $\cal M$.

Given such a triple, we can recover the original pair of ordinary factorization systems by defining $\widetilde{\cal M}$ to be the compositions of morphisms in $\cal F$ followed by morphisms in $\cal M$, and likewise $\widetilde{\cal E}$ as compositions of morphisms in $\cal E$ and $\cal F$.  We will switch between both of these views of a ternary factorization system as convenient.

The 2-category $\cat$ of small categories and functors has a well-known 
ternary factorization system $(\cal E, \cal F, \cal M)$ with $\cal E$ 
consisting of all essentially surjective and full functors, $\cal F$ 
consisting of all essentially surjective and faithful functors, and 
$\cal M$ consisting of all fully-faithful functors. This factorization 
system restricts to small groupoids as well. That is to say, given a 
functor $F: \cal G \to \cal H$ between small groupoids, we can factor 
$F$ in $\cat$ as
\[
G \xrightarrow{F_2} \im_2 F \xrightarrow{F_1} \im_1 F \xrightarrow{F_0} \cal H
\]
and the intermediate categories $\im_2 F$ and  $\im_1 F$ will be groupoids. Since it will be useful later, we will give an explicit construction of this factorization and the intermediate categories for an arbitrary functor $F: \cal G \to \cal H$ between categories.

The 2-coimage of $F$ is the functor $F_2: \cal G \to \im_2 F$. The objects of $\im_2 F$ are the objects of $\cal G$ and the morphisms are equivalence classes of morphisms in $\cal G$ where $f \sim f'$ if and only if $F(f) = F(f')$. The functor $F_2$ acts as the identity on objects and sends morphisms to their corresponding equivalence class.

The 1-image of $F$ is the functor $F_0: \im_1 F \to \cal H$, where 
$\im_1 F$ is the full subcategory of $\cal H$ consisting of all objects 
$Y$ such that $Y \cong F(X)$ for some object $X$ of $\cal G$. The 
functor $F_0$ is the inclusion of this subcategory. 

The functor $F_1: \im_2 F \to \im_1 F$ sends each object $X$ to $F(X)$ and each equivalence class of morphisms $[f]$ to $F(f)$. 

\section{Cardinalities of functor groupoids} \label{sec_funct_grpd_card}
Let $\cal H$ be a finite groupoid and $\cal G$ a locally finite 
groupoid. Then the functor groupoid $\mathcal{G^H}$ is locally finite 
and
\begin{gather}
    |\mathcal{G^H}| = \prod_{[y] \in \overline{\cal H}} \sum_{[x] \in 
    \overline{\cal G}}  \frac{\#\Hom(\cal H_y, \cal G_x)}{\#\cal G_x}
\end{gather}
This formula is not very useful for computing cardinalities, as it 
requires counting homomorphisms between arbitrary finite groups, which 
is a hard problem in general. However, if $\cal G$ is the groupoid of 
finite sets and bijections, then Corollary \ref{cor_g_set_card} gives a 
nice formula for the cardinality. In the particular case where $\cal H 
\simeq BG$ for some finite group $G$, this groupoid is equivalent to 
the core of the category of finite $G$-Sets, which comes with a 
forgetful functor $U: \core(G-\finset) \to \finset$, making 
it a stuff type. We can then 
provide a nice formula for the generating function of $U$. Then the 
cardinality of the functor groupoid $\finset^{BG}_{\cong}$ is 
given by evaluating this generating function at 1.

\theorem
    Let $G$ be a finite group. Then the generating function of the 
    stuff type $U: \core(G\mbox{-}\finset) \to \finset$ is
    \[
        |U|(x) = \exp\left(\sum_{[H\leq 
        G]}\frac{x^{[H:G]}}{\#C_{\sym(G/H)}(\im\,\theta)}\right),
    \]
    where the sum is taken over subgroups $H \leq G$ up to conjugacy, 
    $\sym(G/H)$ is the group of permutations of the set $G/H$, and 
    $\theta: G \to \sym(G/H)$ is the homomorphism defining the usual 
    action of $G$ on left cosets.
\endtheorem
\pf
    First, for each $H \leq G$ we want to find the generating function 
    $\Phi_H$ of finite $G$-sets where all orbits are isomorphic to 
    $G/H$. Clearly, the cardinality of such $G$-sets must be a multiple 
    of $[H:G]$, so the coefficients of $\Phi_H$ are supported only on 
    $x^{[H:G]n}$. In order to compute these coefficients, we need to 
    find the order of $\aut(G/H \times [n])$. Each such automorphism is 
    described by the data of a permutation $\pi$ of the $n$ orbits, as 
    well as an $n$-tuple $(\phi_i)$ of automorphisms of $G/H$. Since 
    the automorphism group of $G/H$ is $C_{\sym(G/H)}(\im\,\theta)$, we 
    have
    \[
    \#\aut(G/H \times [n]) = n!(\#C_{\sym(G/H)}(\im\,\theta))^n.
    \]
    Therefore, the formula for the generating function $\Phi_H$ is
    	$$\mld \Phi_H(x) &= \sum_{n=0}^\infty 
    	\frac{x^{[H:G]n}}{n!(\#C_{\sym(G/H)}(\im\,\theta))^n}\\ = 
        \exp\left(\frac{x^{[H:G]}}{\#C_{\sym(G/H)}(\im\,\theta)}\right).
        $$
    Now, in order to put a $G$-action on an arbitrary finite set, we can partition it into multiple disjoint (possibly empty) subsets - one for each conjugacy class $[H]$ of subgroups of $G$ - and give each one the structure of a $G$-set where all orbits are isomorphic to $G/H$. This operation corresponds to taking the product of generating functions, so we have:
    	$$\mld |U|(x) &= \prod_{[H\leq G]} \Phi_H(x)\\
    	= \exp\left(\sum_{[H\leq 
    	G]}\frac{x^{[H:G]}}{\#C_{\sym(G/H)}(\im\,\theta)}\right).$$
\epf

\corollary \label{cor_g_set_card}
    For a finite groupoid $\cal G$, the cardinality of $\finset^{\cal G}_{\cong}$ is
    \[
    \exp\left(\sum_{[x]\in \overline{\cal G}}\sum_{[H\leq\cal G_x]} \frac1{\#C_{\sym(\cal G_x/H)}(\im\,\theta)}\right)
    \]
\endcorollary

This result uses the so called convolution product or Cauchy product of 
stuff types, which relies only on the fact that $\finset$ has 
coproducts. Replacing $\finset$ with another locally finite category 
with coproducts, such as the category $\mathrm{Vect_{f.d.}}(\Fq)$ of 
finite dimensional vector spaces over the field $\Fq$ 
lets us prove a similar result\\

Let $G$ be a finite group and $\Fq$ a finite field whose 
characteristic does not divide the order of $G$. The groupoid of all 
finite representations of $G$ over $\Fq$ has a functor $U : 
\core(\mathrm{Rep_{f.d.}}(G,\Fq)) \to 
\mathrm{Vect_{f.d.}}(\Fq)$ sending each representation to it's 
underlying vector space. We can define the generating function of $U$ 
similarly by
\[
|U|(x) = \sum_{n=0}^\infty |\cal V_n|x^n,
\]
where $\cal V_n$ is the subgroupoid of all representations of dimension 
$n$.

\theorem
Let $G$ a finite group, and $q = p^k$ for $p$ a prime not dividing the 
order of $G$. Then the groupoid 
$\core(\mathrm{Rep_{f.d.}}(G,\Fq))$ is tame.
\endtheorem
\pf
It is enough to show that the generating function $|U|(x)$ converges at 
$x=1$. Maschke's theorem states that any representation of $G$ is 
isomorphic to a direct sum of irreducible representations, of which 
there are only finitely many. So the generating function $|U|(x)$ can 
be written as
\[
	|U|(x) = \prod_{V} \Phi_V(x),
\]
where the product is taken over irreducible representations of $G$ over 
$\Fq$ up to isomorphism, and $\Phi_V(x)$ is the generating 
function for all representations whose only components are the 
irreducible representation $V$. $\Phi_V$ is given by the formula:
\[
	\Phi_V(x) = \sum_{n=0}^\infty \frac{x^{n \dim 
	V}}{\#\aut(V^{\oplus n})}
\]

Since there are only finitely many irreducible representations $V$ up 
to isomorphism, it is enough to show that each $\Phi_V(x)$ converges at 
$x=1$. To do this, we will find an upper bound for $\Phi_V(1)$. 

Choosing a basis for $V$, we can consider $\aut(V)$ to consist 
of invertible matrices with respect to this basis. Then we can consider
$\aut(V^{\oplus n})$ to be invertible block matrices where the 
blocks are either zero or an element of $\aut(V)$. For a 
fixed $V$, let $B_n \leq \aut(V^{\oplus n})$ be the subgroup of 
block upper triangular matrices. Let $\#\aut(V) = a$, then for 
$n \geq 2$, 
\[
\#B_n = a^n(a+1)^{T_{n-1}}
\]
where $T_n$ are the triangular numbers. Then we can estimate 
$\Phi_V(1)$ by
	$$\mld \Phi_V(1) &= \sum_{n=0}^\infty \frac1{\#\aut(V^{\oplus 
	n})}\\ \leq \sum_{n=0}^\infty \frac1{\#B_n}\\ = 1 + \frac1a + 
	\sum_{n=2}^\infty \frac1{a^n(a+1)^{T_{n-1}}}$$
Since the terms of this series are bounded above by 
$\frac1{(a+1)^{n^2}}$ asymptotically, it must converge. Therefore 
$\Phi_V(x)$ converges at 
$x=1$ for all $V$, and so $\core(\mathrm{Rep_{f.d.}}(G,\Fq))$ is tame.
\epf

When the characteristic of $\Fq$ divides the order of $G$, 
the group algebra decomposes as a sum over the blocks of $G$, rather 
than the irreducible representations. So the generating function will 
be a product of the generating functions for each block. Showing that 
these generating functions converge is much more complicated than the 
case where $p \nmid \#G$. However, it seems likely that a similar 
strategy of finding a lower bound for the growth of the automorphism 
groups that is at least geometric could still work, leading to the 
following conjecture.

\begin{conj}
For any finite group $G$ and finite field $k$, the groupoid 
$\core(\mathrm{Rep_{f.d.}}(G,k))$ is tame.
\end{conj}

\section{Morphisms between tame groupoids} \label{sec_morphisms_tame}
In this section, we prove some facts relating properties of functors to 
the groupoid cardinalities of their source and target groupoids, 
analogous to the relationship between properties of functions and 
cardinality of sets.
\proposition \label{tame_grpd_order}
    Let $\varphi: \cal G \to \cal H$ be a functor between tame groupoids.
    \begin{enumerate}
        \item[(a)] If $\varphi$ is full, then $|\cal G| \leq |\cal H|$.
        \item[(b)] If $\varphi$ is essentially surjective and faithful, 
        then $|\cal G| \geq |\cal H|$.
    \end{enumerate}
\endproposition

\pf
    \textbf{(a)} Since $\varphi$ is full, it induces an injective map on isomorphism classes. So,  the formula for $|\cal H|$ gives
    \[
    |\cal H| = \sum_{[x] \in \overline{\cal G}} \frac1{\#\cal 
    H_{\varphi(x)}} + \sum_{[y] \in \overline{\cal H} \smallsetminus 
    \varphi(\overline{\cal G})} \frac1{\#\cal H_y}.
    \]
    In addition, $\varphi$ induces surjective maps $\cal G_x \to \cal H_{\varphi(x)}$, so
    \[
    \frac1{\#\cal G_x} \leq \frac1{\#\cal H_{\varphi(x)}}
    \]
    Therefore,
    \[
    |\cal G|  = \sum_{[x] \in \overline{\cal G}} \frac1{\#\cal G_x} 
    \leq \sum_{[x] \in \overline{\cal G}} \frac1{\#\cal H_{\varphi(x)}} 
    = |\cal H|.
    \]

    \textbf{(b)} Since $\varphi$ is essentially surjective, it induces a surjective map on isomorphism classes, so the formula for $|\cal G|$ gives
    \[
    |\cal G|  = \sum_{[y] \in \overline{\cal H}} \left(\sum_{[x]\in\varphi^{-1}([y])} \frac 1{\#\cal G_x}
    \right)
    \]
    Since $\varphi$ is faithful, for each $[x]$ in $\overline{\cal G}$, $\varphi$ induces an injective group homomorphism $\cal G_x \to \cal H_{\varphi(x)}$, so
    \[
    \frac1{\#\cal G_x} \geq \frac1{\#\cal H_{\varphi(x)}}
    \]
    Therefore,
    \[
    |\cal G|  = \sum_{[y] \in \overline{\cal H}} \left(\sum_{[x]\in\varphi^{-1}([y])} \frac 1{\#\cal G_x}\right) \geq  \sum_{[y] \in \overline{\cal H}} \frac{\#\varphi^{-1}([y])}{\#\cal H_y} \geq  \sum_{[y] \in \overline{\cal H}} \frac 1{\#\cal H_y} = |\cal H|.
    \]
\epf

\theorem \label{equi_if_equal_card}
Let $\cal G, \cal H$ be tame groupoids such that $|\cal G| = |\cal H|$, and let $\varphi: \cal G \to \cal H$ be a functor that is either
\begin{enumerate}
    \item[(a)] essentially surjective and full,
    \item[(b)] essentially surjective and faithful, or
    \item[(c)] fully faithful.
\end{enumerate}
Then $\varphi$ is an equivalence.
\endtheorem
\pf
\begin{enumerate}
    \item [(a)] Since $\varphi$ is essentially surjective and full, it induces a bijection on isomorphism classes. So since $|\cal G| = |\cal H|$,
    \[
    |\cal G| = \sum_{[x] \in \overline{\cal G}} \frac 1{\#\cal G_x} = \sum_{[x] \in \overline{\cal G}} \frac 1{\#\cal H_{\varphi(x)}} = |\cal H|
    \]
    Additionally, since $\varphi$ is full, it induces surjective group homomorphisms on all vertex groups, which means for all objects $x$ of $\cal G$,
    \[
    \frac 1 {\#\cal G_x} \leq \frac 1 {\#\cal H_{\varphi(x)}}
    \]
    Suppose for contradiction that this inequality is strict for some object $x_0$ of $\cal G$. So
    \[
    0 < \frac 1{\#\cal H_{\varphi(x_0)}} - \frac1{\#\cal G_{x_0}} = \sum_{\substack{[x]\in \overline{\cal G}\\ x \ncong x_0}} \left(\frac1{\#\cal G_x} - \frac1{\#\cal H_{\varphi(x_0)}}\right)
    \]
    However, this is a contradiction, since every term in the sum on 
    the right must be non-positive. Therefore the induced map on each 
    vertex group must be a bijection, making $\varphi$ an equivalence.
    \item [(b)] Since $\varphi$ is essentially surjective, it induces a surjection on isomorphism classes. So,
    \[
    |\cal G| = \sum_{[y]\in\overline{\cal H}}\sum_{\substack{[x] \in \overline{\cal G}\\ \varphi(x) \cong y}} \frac 1{\#\cal G_x} = \sum_{[y]\in\overline{\cal H}} \frac 1{\#\cal H_y} = |\cal H|.
    \]
    Since $\varphi$ is faithful, the induced maps on vertex groups are injective, so for $\varphi(x)\cong y$,
    \begin{gather*}
        \frac1{\#\cal H_y} \leq \frac1{\#\cal G_x}\, \text{, and}\\
        \frac1{\#\cal H_y} \leq \sum_{\substack{[x] \in \overline{\cal G}\\ \varphi(x) \cong y}} \frac 1{\#\cal G_x}.
    \end{gather*}
    Using the same argument as part (a), we see that the second inequality cannot be strict. Following from the first inequality, the sum on the right must contain only one term, which is equal to the left hand side. Thus the induced maps on vertex groups are all isomorphisms, and so $\varphi$ is full, making it an equivalence.
    \item[(c)] Since $\varphi$ is fully faithful, we can treat $\cal G$ as a full subgroupoid of $\cal H$, thus
    \[
    |\cal G| = \sum_{[x]\in\overline{\cal G}}\frac1{\#\cal G_x} = \sum_{[x]\in\overline{\cal G}}\frac1{\#\cal G_x} + \sum_{[y] \in \overline{\cal H}\smallsetminus\overline{\cal G}} \frac1{\#\cal H_y} = |\cal H|
    \]
    Clearly, for these sums to be equal, $\overline{\cal 
    H}\smallsetminus\overline{\cal G} = \varnothing$, so $\varphi$ is 
    essentially surjective, and therefore an equivalence.
\end{enumerate}
\epf

\theorem \label{equi_if_tame}
    Let $\mathcal {G,H}$ be tame groupoids, and let $\varphi: \cal G \to \cal H$ and $\psi: \cal H \to \cal G$ be functors.
    \begin{enumerate}
        \item[(a)] If $\varphi$ and $\psi$ are both full, then they are both equivalences.
        \item[(b)] If $\varphi$ and $\psi$ are both essentially surjective and faithful, then they are both equivalences.
    \end{enumerate}
\endtheorem

\pf
    \textbf{(a)} We want to show that $\psi\phi$ is essentially 
    surjective, so fix an object $x_0$ of $\cal G$. We will find an 
    object $y$ such that $\psi\phi(y)\cong x_0$. Let $x_n = 
    (\psi\varphi)^n(x_0)$. Since $\psi\varphi$ is full, we get a chain 
    of surjective group homomorphisms on the vertex groups of the $x_n$:
    \[
    \cal G_{x_0} \twoheadrightarrow \cal G_{x_1} \twoheadrightarrow \cal G_{x_2} \twoheadrightarrow \dots
    \]
    Since all of these groups are finite, the chain eventually 
    stabilizes. So for some $M$, and all $m \geq M$, $\cal G_{x_m} 
    \cong \cal G_{x_{m+1}}$. If we assume that no two of the $x_n$ are 
    isomorphic in $\cal G$, then $\cal G$ contains an infinite number 
    of components with vertex groups all the same size, which 
    contradicts the assumption that $\cal G$ is tame. Thus, it must be 
    the case that there is some $n$ and $m$ such that $x_{m+n} \cong 
    x_m$. Then, since full functors induce injective functions on 
    isomorphism classes, and
    \[
    (\psi\phi)^m(x_n) = x_{m+n} \cong x_m = (\psi\phi)^m(x_0)
    \]
    we conclude that $x_n \cong x_0$. This implies that $\psi\varphi$ 
    is essentially surjective, and it then follows that $\psi$ is 
    essentially surjective. Also, since $\cal G_{x_0}$ and $\cal 
    G_{x_n}$ are isomorphic finite groups, it must be that each $\cal 
    G_{x_k} \to \cal G_{x_{k+1}}$ is actually a bijection, so 
    $\psi\varphi$ is faithful, and it then follows that $\varphi$ is 
    faithful. We can use the same procedure to show that $\varphi$ is 
    essentially surjective and $\psi$ is faithful. Thus $\varphi$ and 
    $\psi$ are full, faithful, and essentially surjective, and 
    therefore are equivalences.

    \textbf{(b)} Proposition~\ref{tame_grpd_order} says that $|\cal G| 
    \geq |\cal H|$ and $|\cal H| \geq |\cal G|$. Therefore $|\cal G| = 
    |\cal H|$. Then Theorem~\ref{equi_if_equal_card} states that 
    $\varphi$ and $\psi$ are equivalences.
\epf

\section{Relatively finite functors} \label{sec_rel_fin}
A large class of (2,1)-categories with nice combinatorial properties arise as subcategories of the slice (2,1)-category $\grpd/\cal B$ consisting of functors having a certain finiteness property. Before we discuss this property, it's worth discussing the category $\grpd/\cal B$ in more detail.

The objects of $\grpd/\cal B$ are pairs $(\cal G, F)$, where $\cal G$ 
is a groupoid, and $F:\cal G \to \cal B$ is a functor. The 1-morphisms 
from $(\cal G_1, F_1)$ to $(\cal G_2, F_2)$ are pairs $(\Phi,\eta)$, 
with $\Phi$ a functor and $\eta$ a natural isomorphism forming a 
diagram in $\grpd$:
\[\begin{tikzcd}
	{\cal G_1} && {\cal G_2} \\
	\\
	& {\cal B}
	\arrow["\Phi", from=1-1, to=1-3]
	\arrow[""{name=0, anchor=center, inner sep=0}, "{F_1}"', from=1-1, to=3-2]
	\arrow["{F_2}", from=1-3, to=3-2]
	\arrow["\eta"'{pos=0.4}, shorten <=9pt, shorten >=14pt, Rightarrow, from=0, to=1-3]
\end{tikzcd}\]
Given 1-morphisms $(\Phi_1,\eta_1)$ and $(\Phi_2,\eta_2)$ from $(\cal G, J)$ to $(\cal H, K)$ are natural isomorphisms $\psi: \Phi_1 \Rightarrow \Phi_2$ such that $\eta_2 = K\phi \circ \eta_1$.

\remark
Since $\cal B$ is a groupoid, the natural transformation in the 
definition of the 1-morphisms is necessarily an isomorphism. However, 
loosening this requirement (so that $\cal B$ is a general category) 
does not change this definition. The slice 2-category 
$\cat/\cal B$ is also defined so that the natural 
transformations occurring in the definition of 1-morphisms be 
isomorphisms. If we allow $\eta : F_1 \Rightarrow F_2\Phi$ to not be an 
isomorphism, we obtain what is known as a lax morphism in 
$\cat/\cal B$. There are interesting examples of lax morphisms 
arising in combinatorial contexts, for instance, when taking the 
derivative of a stuff type, but we will not consider them here.
\endremark

\definition
    Let $F: \cal G \to \cal B$ be a functor. For any object $y$ of $\cal B$, the \textbf{full inverse image} $F^{-1}(y)$ of $F$ over $y$ is the full subgroupoid of $\cal G$ consisting of all objects $x$ such that $F(x) \cong y$.
    
    If we restrict $F$ to $F^{-1}(y)$, we get a functor $F|_y : F^{-1}(y) \to \cal B$ called the \textbf{component of} $F$ \textbf{at} $x$.
\enddefinition
It is clear that any such $F$ can be written uniquely as a coproduct:
\[
F \cong \coprod_{[x] \in \overline{\cal B}} F|_x
\]
Furthermore, given any 1-morphism $(\Phi,\eta): (F_1,\cal G_1) \to (F_2,\cal G_2)$ in $\grpd/\cal B$ and any object $x$ in $\cal B$, it is clear that $\Phi$ must send $F_1^{-1}(x)$ to $F_2^{-1}(x)$. So restricting $\Phi$ to the components of $F_1$ gives a morphism between components $(\Phi|_x,\eta|_x) : F_1|_x \to F_2|_x$, and therefore $\Phi$ can be written as
\[
\Phi \cong \coprod_{[x]\in \overline{\cal B}} \Phi|_x
\]

\definition
    The functor $F$ is \textbf{relatively finite} if for every object $y$ of $\cal G$, $F^{-1}(y)$ is equivalent to a finite groupoid.
\enddefinition
\definition
    The (2,1)-category $\relfin_{\cal B}$ is the full subcategory of $\grpd/\cal B$ consisting of objects $(\cal G, F)$ such that $F$ is relatively finite.
\enddefinition

\example
    Let $*$ be the discrete groupoid with a single object. Then $\relfin_*$ is equivalent to the category of finite groupoids.
\endexample

\example
    Let $\cal B$ be the category with objects all finite sets and 
    morphisms all bijections between them. Then $\relfin_{\cal B}$ is 
    the category of relatively finite stuff types. Most of the stuff 
    types that arise in combinatorial contexts are relatively finite.
\endexample

\remark
$\relfin_{\cal B}$ has all finite coproducts and if $\cal B$ is locally 
finite, all finite products. In addition if $\cal B$ is the groupoid of 
finite sets and bijections, then $\relfin_{\cal B}$ is closed under the 
usual operations of stuff-types (i.e. the Cauchy product and 
derivative).
\endremark

\subsection{Factorization of morphisms in $\grpd/\cal B$} \label{sec_fact_in_grpd_b}
Since a morphism in $\grpd/\cal B$ is defined as a pair $(F,\eta)$ with $F$ a functor, we can factorize $F$ as just a functor of groupoids. The question is: does this factorization give a factorization of the morphism $(F,\eta)$?

\proposition
Let $(F,\eta):(\cal G, \Phi) \to (\cal H, \Psi)$ be a 1-morphism in 
$\grpd/\cal B$, with $F$ factorizing as $F = F_0\circ F_1\circ F_2$ as 
described in Section~\ref{factorization}. Then there exist functors and 
natural transformations forming the following diagram in $\grpd$:
\[\begin{tikzcd}
	{\cal G} && {\im_2 F} && {\im_1 F} && {\cal H} \\
	\\
	\\
	&&& {\cal B}
	\arrow["{F_2}", from=1-1, to=1-3]
	\arrow[""{name=0, anchor=center, inner sep=0}, "\Phi"', from=1-1, to=4-4]
	\arrow["{F_1}", from=1-3, to=1-5]
	\arrow[""{name=1, anchor=center, inner sep=0}, "{\hat\Phi}"', from=1-3, to=4-4]
	\arrow["{F_0}", from=1-5, to=1-7]
	\arrow[""{name=2, anchor=center, inner sep=0}, "{\hat\Psi}"'{pos=0.6}, from=1-5, to=4-4]
	\arrow["\Psi", from=1-7, to=4-4]
	\arrow["{\eta_2}", shorten <=9pt, shorten >=6pt, Rightarrow, from=0, to=1-3]
	\arrow["{\eta_1}", shorten <=15pt, shorten >=15pt, Rightarrow, from=1, to=1-5]
	\arrow["{\eta_0}", shorten <=23pt, shorten >=31pt, Rightarrow, from=2, to=1-7]
\end{tikzcd}\]
\endproposition
\pf
First, to define $\hat\Phi$, recall that the objects of $\im_2 F$ are 
the objects of $\cal G$, so $\hat\Phi(x) := \Phi(x)$. For an 
equivalence class of morphisms $[f]$, define $\hat\Phi([f]) := 
\Phi(f)$. To see that this is well defined, suppose we have $F(f) = 
F(f')$ for some distinct $f,f':x\to y$. Since $\eta$ is a natural 
isomorphism, we have
	$$\mld
    \Phi(f) &= \eta_y^{-1}\circ \Psi F(f) \circ \eta_x\\
	= \eta_y^{-1}\circ \Psi F(f') \circ \eta_x\\
    = \Phi(f')$$
Then the components of $\eta_2$ are simply the identity morphisms of $\Phi(x)$ for each object $x$ of $\cal G$.

Since $\im_1 F$ is considered as a full subcategory of $\cal H$, simply define $\hat\Psi$ to be the restriction of $\Psi$ to this subcategory. Then the components of $\eta_0$ are the identity morphisms of $\Psi(y)$ for each object $y$ of $\im_1 F$.

Finally, the component of $\eta_1$ corresponding to $x \in \mathrm{Ob}(\im_2 F) = \mathrm{Ob}(\cal G)$, must go from $\hat\Phi(x)= \Phi(x)$ to $\hat\Psi F_1(x) = \Psi F(x)$. The obvious choice is to take this component to be $\eta_x$. Then since $\hat\Phi([f]) = \Phi(f)$ and $\hat\Psi F_1([f]) = \Psi F(f)$, the fact that $\eta$ is natural ensures that $\eta_1$ is natural.
\epf

It remains to show that $(\cal E, \widetilde{\cal M})$ and $(\widetilde{\cal E}, \cal M)$ are factorization systems of $\grpd/\cal B$ when extended this way. The first three axioms are obvious. The orthogonality axiom guarantees the existence of a functor acting as a fill-in, but we must still construct a natural transformation making that functor a morphism in $\grpd/\cal B$. The construction of this natural transformation is again a straightforward, yet tedious, use of the fill-in data from $\grpd$, and is left as an exercise to the reader.

So the standard ternary factorization system on groupoids extends nicely to $\grpd/\cal B$. Moreover, the nice properties of Definition~\ref{fully_co_faith} are also preserved.

\proposition \label{prop_nice_fact}
Let $(G,\varphi): (X_1,F_1) \to (X_2,F_2)$ be a 1-morphism in $\grpd/\cal B$.
\begin{enumerate}
    \item[(a)] If $G$ is fully faithful as a functor, then $(G,\varphi)$ is fully faithful in the sense of Definition~\ref{fully_co_faith}.
    \item[(b)] If $G$ is essentially surjective and full, then $(G,\varphi)$ is fully cofaithful.
\end{enumerate}
\endproposition
\pf
\begin{enumerate}
    \item[(a)] Fix an object $(Y,H)$ in $\grpd/\cal B$. The induced functor is given by
    \begin{align*}
        (G,\varphi)_*: \grpd/\cal B((Y,H),(X_1,F_1)) &\to \grpd/\cal 
        B((Y,H),(X_2,F_2))\\
        (\Phi,\eta) &\mapsto (G\Phi, \varphi\Phi\circ\eta)\\
        (\nu:(\Phi_1,\eta_1) \Rightarrow (\Phi_2,\eta_2)) &\mapsto (G\nu:(G\Phi_1,\varphi\Phi_1\circ\eta_1) \Rightarrow (G\Phi_2,\varphi\Phi_2\circ\eta_2))
    \end{align*}
    Fix $(\Phi_1,\eta_1)$ and $(\Phi_2,\eta_2)$ from $(Y,H)$ to 
    $(X_1,F_1)$, which are sent to $(G\Phi_1, \varphi\Phi_1\circ 
    \eta_1)$ and $(G\Phi_2, \varphi\Phi_2\circ \eta_2)$ by 
    $(G,\varphi)_*$. Since $G$ is fully faithful as a functor,  
    Proposition~\ref{dv_prop_proper_fact} states that $G_*:\grpd(Y,X_1) 
    \to \grpd(Y,X_2)$ is fully faithful as a functor. So given a 
    2-morphism $\mu: (G\Phi_1,\varphi\Phi_1\circ\eta_1) \Rightarrow 
    (G\Phi_2,\varphi\Phi_2\circ\eta_2)$ in $\grpd/\cal B$, there is a 
    unique natural transformation $\nu: \Phi_1 \Rightarrow \Phi_2$ such 
    that $\mu = G\nu$. It remains to show that $\nu$ is a 2-morphism in 
    $\grpd/\cal B$. That is, we want to show that $\eta_2 = F_1\nu 
    \circ \eta_1$.

    Since $\mu$ is a 2-morphism in $\grpd/\cal B$, we have the following identity:
    \[
    \varphi\Phi_2\circ\eta_2 = F_2\mu \circ \varphi\Phi_1 \circ \eta_1
    \]
    Since all 2-morphisms are invertible, using cancellation, it is enough to show that
    \[
    \varphi\Phi_2 \circ F_1\nu = F_2G\nu \circ \varphi\Phi_1
    \]
    We will do this by comparing the components at each object in $Y$. So, fix an object $a$ of $Y$, and consider the morphism $\nu_a$. Since $\varphi: F_1 \Rightarrow F_2G$ is a natural transformation, we have the identity
    \[
   	\varphi_{\Phi_2(a)} \circ F_1(\nu_a) = F_2G(\nu_a) \circ 
   	\varphi_{\Phi_1(a)}
    \]
    Then by definition of the whiskerings $\varphi\Phi_1$, $\varphi\Phi_2$, $F_1\nu$, and $F_2G\nu$, we have shown that $\nu$ is a 2-morphism in $\grpd/\cal B$ from $(Y,H)$ to $(X_1,F_1)$.

    \item[(b)] Fix an object $(Y,H)$ in $\grpd/\cal B$. The induced functor is given by
    \begin{align*}
        (G,\varphi)^*: \grpd/\cal B((X_2,F_2),(Y,H)) &\to \grpd/\cal B((X_1,F_1),(Y,H))\\
        (\Phi,\eta) &\mapsto (\Phi G, \eta G\circ \varphi)\\
        (\nu:(\Phi_1,\eta_1) \Rightarrow (\Phi_2,\eta_2)) &\mapsto (\nu G:(\Phi_1 G,\eta_1 G \circ \varphi) \Rightarrow (\Phi_2 G,\eta_2 G \circ \varphi))
    \end{align*}
    Fix $(\Phi_1,\eta_1)$ and $(\Phi_2,\eta_2)$ from $(X_2,F_2)$ to $(Y,H)$, and let $\mu: (\Phi_1G,\eta_1G \circ \varphi) \Rightarrow (\Phi_2G,\eta_2G\circ\varphi)$ be a 2-morphism in $\grpd/\cal B$, so that
    \[
        H\mu\circ \eta_1G\circ \varphi = \eta_2G\circ \varphi.
    \]
    Then since $\varphi$ is invertible, using right cancellation gives
    \[
        H\mu \circ \eta_1G = \eta_2G
    \]
    Proposition~\ref{dv_prop_proper_fact} states that there is a unique natural transformation $\nu:\Phi_1 \Rightarrow \Phi_2$ such that $\mu = \nu G$. Then,
    \[
    (H\nu \circ \eta_1)G = \eta_2 G
    \]
    Then using Proposition~\ref{dv_prop_proper_fact} again, we see that
    \[
    H\nu \circ \eta_1 = \eta_2
    \]
    Which proves that $\nu$ is a 2-morphism in $\grpd/\cal B$.
\end{enumerate}
\epf

\section{Combinatorial (2,1)-categories} \label{sec_comb_cat}
In this section we prove the main theorem of this paper, which is a generalization of Theorem \ref{lovasz_lemma} to categories of relatively finite functors. The first half of the proof is similar to the proof of the original theorem. We will need the following definition and lemma:
\definition
    Let $\cal C$ be a (2,1)-category with a factorization system $(\cal 
    E,\cal M)$. For any object $X$ of $\cal C$, the $\cal 
    E$-\textbf{quotients} of $X$ are equivalence classes of morphisms 
    $p: X \to Y$ in $\cal E$, where two morphisms $p:X \to Y$ and $p':X 
    \to Y'$ are considered equivalent if there exists an equivalence 
    $f: Y \to Y'$ and a 2-morphism $\varphi: fp \Rightarrow p'$.
\enddefinition

\lemma \label{lem_factor_decomp}
    Let $\cal C$ be a (2,1)-category with a factorization system $(\cal E, \cal M)$ such that all 1-morphisms in $\cal E$ are fully cofaithful. Then for any two objects $X,Y$, there is an equivalence of groupoids
    \[
    \cal C(X,Y) \simeq \coprod_{[X\twoheadrightarrow Z]\in\cal E}\cal M(Z,Y)
    \]
    where the coproduct is taken over all $\cal E$-quotients of $X$.
\endlemma
\pf
    Choosing representatives for the $\cal E$-quotients of $X$, we can construct a functor from the coproduct on the right to $\cal C(X,Y)$ by pre-composing with the corresponding representatives. This functor is fully faithful when restricted to each component of the coproduct, and thus on the entire coproduct. It is obviously essentially surjective as a result of the factorization system.
\epf

This lemma is essential for the first part of the proof of the main 
theorem, as it gives us a formula for $\cal C(X,Y)$ in terms of a 
coproduct of $\cal M(Z,Y)$, which we can convert to a sum of groupoid 
cardinalities. The next lemma is used in the second part and will help 
us find equivalences in $\relfin_{BG}$.

\lemma \label{iso_to_eq}
Let $H_1, H_2$ and $G$ be finite groups and $F_1: BH_1 \to BG, F_2: 
BH_2 \to BG$ be objects of $\relfin_{BG}$. Let $(\varphi, g): F_1 \to 
F_2$ be a 1-morphism in $\relfin_{BG}$ such that $\varphi$ is an 
isomorphism of groups. Then $(\varphi, g)$ is an equivalence in 
$\relfin_{BG}$.
\endlemma
\pf
Treat $F_1:H_1 \to G$ and $F_2:H_2 \to G$ as group homomorphisms. Then 
since $(\varphi, g)$ is a 1-morphism in $\relfin_{BG}$, we can treat 
$g$ as an element of $G$ such that for all $h \in H_1$, $F_1(h) = 
g^{-1}F_2(\varphi(h))g$. Then since $\varphi$ is an isomorphism, for 
every $k \in H_2$, $F_2(k) = gF_1(\varphi^{-1}(k))g^{-1}$. Therefore 
$(\varphi^{-1}, g^{-1})$ is a strict inverse to $(\varphi,g)$.
\epf

\theorem \label{main_thm}
Let $\cal B$ be a locally finite groupoid and let $F: \cal G \to \cal 
B$ and $F':\cal G' \to \cal B$ be functors in $\relfin_{\cal B}$ such 
that for all finite groups $H$ and functors $S: BH \to \cal B$,
\[
|\relfin_{\cal B}(S,F)| = |\relfin_{\cal B}(S,F')|.
\]
Then $F$ and $F'$ are equivalent.
\endtheorem
\pf
First, note that any $F$ in $\relfin_{\cal B}$ is completely determined up to equivalence by the components $F|_y$ for each isomorphism class $[y]$ of $\cal B$. Furthermore, if we have $S: BH \to \cal B$, then any morphism $(S,BH) \to (F,\cal G)$ is completely determined by a morphism $(S,BH) \to (F|_{S*},F^{-1}(S*))$, where $*$ is the unique object of $BH$. Therefore, it is enough to prove this theorem in the case where $\cal B = BG$ for some finite group $G$.\\

The proof has two steps: First, we use induction to show that for all $S:BH \to BG$,
\begin{equation}\label{same_num_faithful}
    |\widetilde{\cal M}(S,F)| = |\widetilde{\cal M}(S,F')|,
\end{equation}
where $\widetilde{\cal M}$ consists of all morphisms whose underlying 
functors are faithful. Second, use equation \ref{same_num_faithful} to 
show that $F$ and $F'$ are of the form $F \simeq U \sqcup V$ and $F' 
\simeq U \sqcup W$ with
\[
|\widetilde{\cal M}(S,V)| = |\widetilde{\cal M}(S,W)|
\]
Since $\cal G$ and $\cal G'$ are finite, we can then use induction to show that $F$ and $F'$ decompose as coproducts with equivalent components, and are therefore themselves equivalent.\\

For the first step, recall that the category $\relfin_{BG}$ has the ternary factorization system $(\cal E, \cal F, \cal M)$ as discussed in Section $\ref{sec_fact_in_grpd_b}$ (equivalently described by the pair of factorization systems $(\cal E,\widetilde{\cal M})$ and $(\widetilde{\cal E},\cal M)$). Furthermore, all morphisms in $\cal E$ are fully-cofaithful. Therefore, we can use Lemma \ref{lem_factor_decomp} with respect to the factorization system $(\cal E, \widetilde{\cal M})$. Applying this lemma to both sides the equation $|\relfin_{BG}(S,F)| = |\relfin_{BG}(S,F')|$ and rearranging terms gives
\begin{equation}\label{ih_sum}
|\widetilde{\cal M}(S,F)| - |\widetilde{\cal M}(S,F')| = \sum_{\substack{[S\twoheadrightarrow T] \in \cal E\\ T \not\simeq S}} |\widetilde{\cal M}(T,F)| - |\widetilde{\cal M}(T,F')|,
\end{equation}

where the sum on the right is taken over $\cal E$-quotients of $S$ that 
are not equivalent to $S$. Note that since $S$ is of the form $BH \to 
BG$, for some finite group $H$, the $\cal E$-quotients of $S$ are of 
the form $BH' \to BG$ for some quotient $H'$ of $H$. So we can define a 
partial order on the equivalence classes of functors $BH \to BG$ for 
$H$ a finite group as follows: Let $S':BH' \to BG \leq S:BH\to BG$ if 
$S'$ is a $\cal E$-quotient of $S$. It is easy to check that this is a 
partial order. Furthermore, since any finite group has only finitely 
many quotients, there can be no infinite descending chains, so this 
partial order is well-founded.\\

Now, fix $S: BH \to BG$, and assume for all quotient groups $H'$ of $H$ 
and functors $T:BH' \to BG$, that $|\widetilde{\cal M}(T,F)| = 
|\widetilde{\cal M}(T,F')|$. Then summands on the right side of 
equation \ref{ih_sum} are all zero, so $|\widetilde{\cal M}(S,F)| = 
|\widetilde{\cal M}(S,F')|$. Then, using the principle of Noetherian 
induction, we can conclude that for any finite group $H$ and functor 
$S:BH \to BG$,
\[
|\widetilde{\cal M}(S,F)| = |\widetilde{\cal M}(S,F')|.
\]
This concludes the first part of the proof.\\

For the second part, assume that $\cal G$ and $\cal G'$ are skeletal, and choose some isomorphism class $BK_0$ of $\cal G$ so that $\cal G \simeq BK_0 \sqcup V_0$ and $F \simeq T_0 \sqcup F|_{V_0}$ for some functor $T_0:BK_0 \to BG$. Clearly, the inclusion of $BK_0$ into $\cal G$ is faithful, and so the morphism $T_0 \to F$ is in $\widetilde{\cal M}$. So
\[
|\widetilde{\cal M}(T_0,F)| = |\widetilde{\cal M}(T_0, F')| \neq 0
\]
and thus there exists a morphism $\varphi_0: T_0 \to F'$ in $\widetilde{\cal M}$. The underlying functor of this morphism must send the unique object of $BK_0$ to some $BJ_0 \subseteq \cal G'$. We can then write $\cal G' \simeq BJ_0 \sqcup W_0$ and $F' \simeq R_0 \sqcup F'|_{W_0}$ for some functor $R_0 : BJ_0 \to BG$. Repeating the same procedure gives a morphism $\psi_0: R_0 \to F$, also in $\widetilde{\cal M}$. Repeating this procedure gives a chain of morphisms in $\widetilde{\cal M}$:
\[
T_0 \xrightarrow{\varphi_0} R_0 \xrightarrow{\psi_0} T_1 \xrightarrow{\varphi_1} R_1 \xrightarrow{\psi_1} T_2 \xrightarrow{\varphi_2} \dots
\]
for $T_i: BK_i \to BG$, $R_i:BJ_i \to BG$ with $F \simeq T_i \sqcup V_i$ and $F' \simeq R_i \sqcup W_i$. Note that we can consider this chain as a chain of injective homomorphisms of finite groups:
\[
K_0 \xrightarrow{\varphi_0} J_0 \xrightarrow{\psi_0} K_1 \xrightarrow{\varphi_1} J_1 \xrightarrow{\psi_1} K_2 \xrightarrow{\varphi_2} \dots
\]
Since $\cal G$ and $\cal G'$ are finite, there are only finitely many 
such $K_i, J_i$ up to isomorphism, so this chain must eventually 
stabilize, and thus there is some $i$ for which $\varphi_i: K_i \to 
J_i$ is an isomorphism, and from Lemma \ref{iso_to_eq} the morphism of 
relatively finite functors $\varphi_i: T_i \to R_i$ is an equivalence. 
Thus we have found a $U:BK \to BG$ such that $F \simeq U \sqcup V$ and 
$F' \simeq U \sqcup W$.\\

Now, for any $S:BH \to BG$, it is easy to see that
\begin{gather*}
    \widetilde{\cal M}(S,F) \simeq \widetilde{\cal M}(S,U) \sqcup \widetilde{\cal M}(S,V)\, \text{and}\\
    \widetilde{\cal M}(S,F') \simeq \widetilde{\cal M}(S,U) \sqcup \widetilde{\cal M}(S,W).
\end{gather*}
So, from equation $\ref{same_num_faithful}$, we see that
	$$\mld |\widetilde{\cal M}(S,V)| &= |\widetilde{\cal M}(S,F)| - 
	|\widetilde{\cal M}(S,U)|\\ = |\widetilde{\cal M}(S,F')| - 
	|\widetilde{\cal M}(S,U)|\\ = |\widetilde{\cal M}(S,W)|$$
We can then use induction to show that $V \simeq W$, and therefore $F \simeq F'$.
\epf

It's worth remarking that the overarching procedure of this proof 
differs slightly from that of Theorem \ref{lovasz_lemma}. The first 
half is the same, where we prove by induction that $|\widetilde{\cal 
M}(S,F)| = |\widetilde{\cal M}(S,F')|$ for all $S$. The second half is 
quite different, however. We can show that there exist 1-morphisms $F 
\to F'$ and $F' \to F$. However, since the underlying functors are 
merely faithful, we cannot use Theorem \ref{equi_if_tame} to show that 
$F \simeq F'$. Instead, we rely on relative finiteness to show that $F$ 
and $F'$ must have connected components that are equivalent, and then 
use induction to show that $F \simeq F'$.

\remark
The general format of this proof can be applied to other 
(2,1)-categories with minor changes. However, finding a general set of 
axioms that characterize all - or even most - combinatorial 
(2,1)-categories seems difficult.
\endremark

\section{Postnikov systems and homotopy cardinality} 
\label{sec_homotopy}
Groupoid cardinality is a special case of the homotopy cardinality of 
an $\infty$-groupoid, for when the $\infty$-groupoid is 1-truncated. 
\definition[\cite{bd}]
Let $X$ be a connected $\infty$-groupoid with all homotopy groups 
finite, then the homotopy cardinality of $X$ is
\[
|X| = \prod_{k=1}^\infty (\#\pi_k(X))^{(-1)^k}
\]
If $X$ is not connected, then the homotopy cardinality is the sum of 
the cardinalities of each connected component.
\enddefinition
It is natural to ask which of the results discussed in this paper can 
be generalized to $\infty$-groupoids.

The prevalence of factorization systems in the above results also leads 
us to the natural generalization of factorization systems to 
$\infty$-categories. Given a morphism of $\infty$-groupoids $f:X \to 
Y$, for each $n \geq 1$, there exists an $\infty$-groupoid $\im_n f$ 
and morphisms $e_n: X \to \im_nf$ and $m_n:\im_nf \to Y$ such that the 
induced maps on homotopy groups satisfy:
\begin{itemize}
    \item $e_n: \pi_{n-1}(X,x) \to \pi_{n-1} (\im_nf,e_n(x))$ is surjective.
    \item $e_n: \pi_{k}(X,x) \to \pi_{k}(\im_nf,e_n(x))$ is a bijection for $k < n-1$.
    \item $m_n: \pi_{n-1}(\im_nf,x) \to \pi_{n-1}(Y,m_n(x))$ is injective.
    \item $m_n: \pi_k(\im_nf,x) \to pi_k(Y,m_n(y))$ is a bijection for $k > n-1$.
\end{itemize}
In other words, $e_n$ is $n$-connected and $m_n$ is $n$-truncated, so this factorization system is referred to as the ($n$-connected, $n$-truncated) factorization system. Furthermore, these classes of morphisms satisfy a higher categorical version of the orthogonality condition in Definition \ref{def_fact_syst}.1 (see Example 5.2.8.16 and Proposition 5.2.8.11 in \cite{htt} ). It is easy to check that if we consider groupoids as $\infty$-groupoids, the (2-connected, 2-truncated) factorization is exactly the (essentially surjective \& full, faithful) factorization, and the (1-connected, 1-truncated) factorization is the (essentially surjective, fully-faithful) factorization system discussed in Section \ref{factorization}.

The orthogonality condition, along with the fact that $n$-connectivity 
implies $(n-1)$-connectivity, gives a morphism $\varphi_n: \im_nf \to 
\im_{n-1}f$ such that $e_{n-1} \simeq \varphi_n \circ e_n$ and $m_n 
\simeq m_{n-1}\circ \varphi_n$. The resulting tower of 
$\infty$-groupoids is called the \emph{relative Postnikov system of} 
$f$ (c.f. \cite[Definition VI.2.9]{goerss-jardine}). This is the 
natural generalization of the ternary factorization system in $\grpd$ 
to $\infty$-groupoids. We thus get a commutative diagram:

\[\begin{tikzcd}[ampersand replacement=\&]
	\&\& X \\
	\dots \& {\im_{n+1}f} \& {\im_nf} \& \dots \& {\im_1 f} \\
	\&\& Y
	\arrow["{e_{n+1}}"', from=1-3, to=2-2]
	\arrow["{e_n}", from=1-3, to=2-3]
	\arrow["{e_1}", from=1-3, to=2-5]
	\arrow[from=2-1, to=2-2]
	\arrow[from=2-2, to=2-3]
	\arrow["{m_{n+1}}"', from=2-2, to=3-3]
	\arrow[from=2-3, to=2-4]
	\arrow["{m_n}", from=2-3, to=3-3]
	\arrow[from=2-4, to=2-5]
	\arrow["{m_1}", from=2-5, to=3-3]
\end{tikzcd}\]

The following result generalizes Proposition \ref{tame_grpd_order}.
\proposition
Let $f : X \to Y$ be a morphism between tame $\infty$-groupoids. Then
\begin{align*}
    |\im_n f| \geq |\im_{n-1} f| &\, \mbox{ if $n$ is even, and}\\
     |\im_n f| \leq |\im_{n-1} f| &\, \mbox{ if $n$ is odd}.
\end{align*}
\endproposition
\pf
In the case where $n = 2$, we have an surjective function $e_1 : \pi_0(X) \twoheadrightarrow\pi_0(\im_1 f)$, and bijections $\pi_k(\im_2 f, y) \cong \pi_k(Y,m_2(y))$ for $k > 1$. So,
	$$\mld |\im_2 f| &= \sum_{x \in \pi_0(X)}\left[\#\pi_1(\im_2 f, 
	e_2(x))^{-1}\cdot \prod_{k=2}^\infty 
	\#\pi_k(Y,f(x))^{(-1)^k}\right]\\ = \sum_{u \in \pi_0(\im_1 
	f)}\left[\sum_{x \in e_1^{-1}(u)}\left(\#\pi_1(\im_2 f, 
	e_2(x))^{-1}\cdot \prod_{k=2}^\infty 
	\#\pi_k(Y,f(x))^{(-1)^k}\right)\right]$$
We have bijections $\pi_k(\im_1 f,y) \cong \pi_k(Y,m_1(y))$ for $k > 0$, so
\[
|\im_1 f| = \sum_{u \in \pi_0(\im_1 f)}\left[\#\pi_1(Y,m_1(u))^{-1}\cdot\prod_{k=2}^\infty \#\pi_k(Y,m_1(u))^{(-1)^k}\right]
\]
So we have a correspondence of the terms in each sum with respect to $\pi_0(\im_1 f)$. When comparing the terms, it is clear that the terms in the formula for $|\im_2 f|$ will have at least as many summands as those for $|\im_1 f|$, which can have at most one term. Furthermore, we have an injective homomorphism $m_2: \pi_1(\im_2(f), e_2(x)) \cong \pi_1(Y, f(x))$, so $\#\pi_1(\im_2 f, e_2(x))^{-1} \geq \#\pi_1(Y,f(x))^{-1}$. Therefore,
\[
|\im_2 f| \geq |\im_1 f|
\]

In the case where $n > 2$, we have $\pi_0(X) \cong \pi_0(\im_n f)\cong \pi_0(\im_{n-1} f)$, and bijections $\pi_k(\im_nf,x)\cong \pi_k(Y,m_n(x))$ for $k > n-1$ and $\pi_k(X,x) \cong \pi_k(\im_nf,e_n(x))$ for $k < n-1$. Therefore,
\begin{multline*}
    |\im_n f| = \sum_{x \in \pi_0(X)}\Bigg[\prod_{k=1}^{n-3} \#\pi_k(X,x)^{(-1)^k}\cdot \#\pi_{n-2}(X,x)^{(-1)^{n-2}}\\ 
    \cdot \#\pi_{n-1}(\im_n f,e_n(x))^{(-1)^{n-1}} \cdot \prod_{k=n}^\infty\#\pi_k(Y,f(x))^{(-1)^k}
    \Bigg]
\end{multline*}
Similarly, in the formula for $|\im_{n-1}f|$, we have bijections $\pi_k(\im_{n-1}f,x) \cong \pi_k(Y,m_{n-1}(x))$ for $k > n-2$ and $\pi_k(X,x)\cong \pi_k(\im_{n-1}f,e_{n-1}(x))$ for $k < n-2$. Therefore,
\begin{multline*}
    |\im_{n-1}f| = \sum_{x \in \pi_0(X)}\Bigg[ \prod_{k=1}^{n-3}\#\pi_k(X,x)^{(-1)^k} \cdot \#\pi_{n-2}(\im_{n-1}f,e_{n-1}(x))^{(-1)^{n-2}} \\ 
    \cdot \#\pi_{n-1}(Y,f(x))^{(-1)^{n-1}} \cdot \prod_{k=n}^\infty\#\pi_k(Y,f(x))^{(-1)^k}\Bigg]
\end{multline*}
In addition, we have a surjective homomorphism $e_{n-1}: \pi_{n-2}(X,x) \twoheadrightarrow \pi_{n-2}(\im_{n-1}f,e_{n-1}(x))$ and an injective homomorphism $m_n: \pi_{n-1}(\im_nf,x) \hookrightarrow \pi_{n-1}(Y,m_n(x))$.\\

Thus, if $n$ is even,
\begin{multline*}
\#\pi_{n-1}(\im_n f,e_n(x))^{(-1)^{n-1}} \cdot \#\pi_{n-2}(X,x)^{(-1)^{n-2}}\\
\geq \#\pi_{n-1}(Y,f(x))^{(-1)^{n-1}} \cdot \#\pi_{n-2}(\im_{n-1}f,e_{n-1}(x))^{(-1)^{n-2}}
\end{multline*}
and therefore $|\im_n f| \geq |\im_{n-1}f|$. Similarly, if $n$ is odd, this inequality is reversed, and we have $|\im_n f| \leq |\im_{n-1}f|$.
\epf

\corollary
If $f: X \to Y$ is a morphism of tame $\infty$-groupoids that is $(n-1)$-connected and $n$-truncated, then
\begin{align*}
    |X| \geq |Y| &\, \mbox{ if $n$ is even}\\
     |X| \leq |Y| &\, \mbox{ if $n$ is odd}.
\end{align*}
\endcorollary
\pf
Since $f$ is $n$-truncated,
\[
X \xrightarrow{\id_X}X \xrightarrow{f} Y
\]
is a valid factorization in the ($n$-connected, $n$-truncated) factorization system, so by orthogonality, $X \simeq \im_nf$. Since $f$ is $(n-1)$-connected,
\[
X \xrightarrow{f} Y \xrightarrow{\id_Y}Y
\]
is a valid factorization in the ($(n-1)$-connected, $(n-1)$-truncated) factorization system, so $\im_{n-1}f \simeq Y$.
\epf

Regarding Theorem \ref{equi_if_tame}, it is clear that if $X$ and $Y$ 
are tame $\infty$-groupoids with only finitely many connected 
components and we have morphisms $X \to Y$ and $Y \to X$ which are 
$(n-1)$-connected and $n$-truncated, then $X$ and $Y$ have all 
isomorphic homotopy groups. However, it remains open whether the 
morphisms in question are indeed equivalences.

The method used to prove Theorem \ref{main_thm} would likely not work when generalizing to $\infty$-groupoids. First of all, the first step relies on the fact that groupoids are already 1-truncated, so we can begin by showing that the cardinality of the groupoids of faithful morphisms into $F$ and $F'$ agree. However, even if we assume that our objects are $n$-truncated, the method used to show that there exist equivalent connected components does not even generalize to 2-groupoids.

\bibliographystyle{plain}
\bibliography{combinatorics_in_2_1_categories}

\end{document}